\documentstyle{amsppt}
\magnification 1200
\def\today
{\ifcase\month\or
     January\or February\or March\or April\or May\or June\or
     July\or August\or September\or October\or November\or December\fi
     \space\number\day, \number\year}
\magnification 1200
\input pictex
\UseAMSsymbols
\hsize 5.5 true in
\vsize 8.5 true in
\parskip=\medskipamount
\NoBlackBoxes

\def\mathbb{\Bbb}

\def\mathcal{\Cal}

\def\supp{\text{\rm supp\,}}

\def\Mat{\text{\rm Mat}}

\def\ve{\varepsilon}
\def\vp{\varphi}
\def\arrowk{^\to{\kern -6pt\topsmash k}}
\def\arrowK{^{^\to}{\kern -9pt\topsmash K}}
\def\arrowr{^\to{\kern-6pt\topsmash r}}
\def\bark{\bar{\kern-0pt\topsmash k}}
\def\arrowvp{^\to{\kern -8pt\topsmash\vp}}
\def\arrowf{^{^\to}{\kern -8pt f}}
\def\arrowg{^{^\to}{\kern -8pt g}}
\def\arrowu{^{^\to}a{\kern-8pt u}}
\def\arrowt{^{^\to}{\kern -6pt t}}
\def\arrowe{^{^\to}{\kern -6pt e}}
\def\tk{\tilde{\kern 1 pt\topsmash k}}
\def\barm{\bar{\kern-.2pt\bar m}}
\def\barN{\bar{\kern-1pt\bar N}}
\def\barA{\, \bar{\kern-3pt \bar A}}

\def\mathbb{\Bbb}

\def\mes{\text {\rm mes\,}}
\TagsOnRight
\NoRunningHeads

\def\supp{\text{\rm supp\,}}

\def\Mat{\text{\rm Mat}}

\def\ve{\varepsilon}
\def\vp{\varphi}
\hsize = 6.2true in

\vsize=8.2 true in
\TagsOnRight
\NoRunningHeads

\document
\topmatter
\title
\smc On the Furstenberg measure and density of states for the Anderson-Bernoulli model at small disorder
\endtitle
\author
J.~Bourgain 
\endauthor
\address
Institute for Advanced Study, Princeton, NJ 08540
\endaddress
\email
bourgain\@ias.edu
\endemail

\abstract
We establish new results on the dimension of the Furstenberg measure and the regularity of the integrated density of states for the Anderson-Bernoulli
model at small disorder.
\endabstract
\endtopmatter

\noindent
{\bf \S0. Summary}

Let $H=\Delta+\lambda V$ where $\Delta$ is the lattice Laplacian on $\Bbb Z$ and $V=(V_n)_{n\in\Bbb Z}$ are independent random variables in
$\{1, -1\}$.
We assume $|\lambda|$ small and restrict the energy $E$ outside a fixed neighborhood of $\{0, 2, -2\}$.
It is shown that the Furstenberg measure $\nu_E$ of the corresponding $SL_2(\Bbb R)$-cocycle $\textstyle\pmatrix E-\lambda V_n &-1\\ 1&0\endpmatrix$ has dimension
at least $\gamma(\lambda)$, where $\gamma(\lambda)\overset {\lambda\to 0}\to \longrightarrow 1$.
As a consequence, we derive that the integrated density of states (IDS) $\Cal N(E)$ is H\"older-regular with exponent at least $s(\lambda)\overset {\lambda\to 0}\to
\longrightarrow 1$.

The spectral theory of the Anderson-Bernoulli (A-B) model has been studied by various authors.
It was shown by Halperin (see [S-T]) that for fixed $\lambda>0$,  $\Cal N(E)$ is not H\"older continuous of any order $\alpha$ larger than
$$
\alpha_0 =\frac{2\log 2}{{\text {Arc \ cosh }} (1+\lambda)}.\tag 0.1
$$
H\"older regularity for some $\alpha>0$ has been established in several papers.
In \hbox{[Ca-K-M]}, le Page's method is used.
Different approaches (including one using the super-symmetric formalism) appear in the important paper [S-V-W] that relies on harmonic analysis
principles around the uncertainty principle.
In [B1], the author proved H\"older regularity of the IDS using the Figotin-Pastur expansion of the Lyapounov exponent and martingale theory.
We note that in both [S-V-W] and [B1], the H\"older exponent $\alpha$ remains uniform for $\lambda\to 0$ (in fact, [B1] gives an explicit exponent
$\alpha(\lambda)> \frac 15+\ve$ for $\lambda\to 0$). 

Thus the result in this Note just falls short of establishing the conjectured Lipschitz regularity of IDS of the A-B model for small $\lambda$.
Related is the question whether the Furstenberg measure on projective space is absolutely continuous when $\lambda$ is small (or even better).
As pointed out at the end of the paper, a natural approach to these problems is through certain spectral gap properties that do not depend on
hyperbolicity.
There have been recent advances (cf. [BG1], [BG2], [B2]), that are based on methods from arithmetic combinatorics.
But presently, this theory seems to restrictive for an application to A-B-cocycles.
It does apply however for Schr\"odinger operators with single site distribution given by a measure of positive dimension.  
\bigskip

\noindent
{\bf \S1 Probabilistic inequalities on the Boolean cube}

The following statement is a consequence of Sperner's combinatorial Lemma.\footnote
"{$^*$}"{It was also used in [B1] and [B-K] in the context of the Anderson-Bernoulli model.}

\proclaim
{Lemma 1} Let $f=f(\ve_1, \ldots, \ve_n)$ be a real valued function on $\{1, -1\}^n$ and denote
$$
I_j= f|_{\ve_j=1} - f|_{\ve_j=-1}\tag 1.1
$$
the $j$-influence, which is a function of $\ve_{j'}, j'\not= j$.

Assume that for all $j=1, \ldots, n$
$$
I_j\geq 0\tag 1.2
$$
(i.e. $f$ is monotone increasing) and moreover
$$
I_j\geq \kappa>0\text { on } \ \Omega_j\cap\Omega_j'\tag 1.3
$$
where $\Omega_j$ (resp. $\Omega_j'$) are subsets of $\{1, -1\}^n$ depending only on the variables $\ve_1, \ldots, \ve_{j-1}$ (resp. $\ve_{j+1}, \ldots, \ve_n$).
Then, for any $t\in \Bbb R$, we have
$$
\mes [|f-t|<\kappa] \leq \frac 1{\sqrt n} +\sum_j (2-\mes \Omega_j -\mes\Omega_j').\tag 1.4
$$
\endproclaim

\noindent
{\it Proof.}

Denote
$$
\tilde\Omega =\bigcap_{1\leq j\leq n} (\Omega_j\cap \Omega_j') \tag 1.5
$$
for which
$$
1-\mes \tilde\Omega \leq \sum_j (2-\mes \Omega_j -\mes\Omega_j').\tag 1.6
$$
We claim that the set $[|f-t|<\kappa] \cap\tilde\Omega$ does not contain a pair of distinct comparable elements $\ve =(\ve_j)_{1\leq j\leq n}$
and $\ve'=(\ve_j')_{1\leq j\leq n}$.
Assume otherwise and $\ve < \ve'$, i.e. $\ve_j\leq \ve_j'$ for each $j$. Then
$$
\align
f(\ve') -f(\ve)&=\sum_{1\leq j\leq n} \big(f(\ve_1, \ldots, \ve_{j-1}, \ve_j',\ldots, \ve_n') -f(\ve_1, \ldots, \ve_j, \ve_{j+1}', \ldots, \ve_n')\big)\\
&=\sum_{\Sb 1\leq j\leq n\\ \ve_j \not= \ve_j'\endSb} I_j(\ve_1, \ldots, \ve_j, \ve_{j+1}', \ldots, \ve_n').\tag 1.7
\endalign
$$
Since  $\ve\in\Omega_j, \ve'\in \Omega_j'$, it follows from our assumption on $\Omega_j, \Omega_{j'}$ that
$$
(\ve_1, \ldots, \ve_j, \ve_{j+1}', \ldots, \ve_n')\in \Omega_j\cap \Omega_j'
$$
and hence $I_j (\ve_{1}, \ldots, \ve_j, \ve_{j+1}', \ldots\ve_n')\geq\kappa$ by (1.3).
In particular, since $\ve \not= \ve'$, 
$$
(1.7) \geq \# \{1\leq j\leq n; \ve_j\not= \ve_j'\}\kappa \geq \kappa
$$
which is however impossible if $|f(\ve) -t|\leq \kappa$ and $|f(\ve')-t|\leq \kappa$.
This establishes the claim.

Therefore, by Sperner's lemma on the maximal size of subsets of $\{1, -1\}^n$  not containing any pair of distinct comparable elements, we get
$$
\mes (\tilde\Omega \cap [|f-t|<\kappa])\lesssim \frac 1{\sqrt n}\tag 1.8
$$
and (1.4) follows from (1.6), (1.8). This proves Lemma 1.

We will use the following corollary of Lemma 1.

\proclaim
{Lemma 2} Let $f$ and $I_j$ be as in Lemma 1 and assume each $I_j\geq 0$.

Assume further $\kappa, \delta>0$ and for each $1\leq j< n$
$$
f|_{ \ve_j=1, \ve_{j+1}=1} -f|_{\ve_j=-1, \ve_{j+1}=-1} \geq \kappa \text { for } \ve\in\Omega_j\tag 1.9
$$
where $\Omega_j\subset\{1, -1\}^n$ is a set only depending on the variables $\ve_{j+2}, \ldots , \ve_n$ and such that
$$
\mes\Omega_j> 1-\delta.\tag 1.10
$$
Then, for all $t\in\Bbb R$
$$
\mes[|f-t|< \kappa]\lesssim \frac 1{\sqrt n}+n\delta.\tag 1.11
$$
\endproclaim

\noindent
{\it Proof.} Assume $n =2 m$ even and write $\omega =(\ve_1, \ve_1', \ldots, \ve_m, \ve_m')$ for the $\{1, -1\}^n$-variable.
With this notation, let $\Omega_j$ refer to the set $\Omega_{2j-1}$.

Partition
$$
\{1, -1\}^{2m} =\bigcup_{S\subset\{1, \ldots, m\}} V_S
$$
with
$$
V_S =\{\omega; \ve_j=\ve_j' \text { if } \ j\in S\text { and } \ \ve_j\not= \ve_j' \text { if } \ j\not\in S\}.\tag 1.12
$$
Thus
$$
\mes [|f- t|<\kappa]=\sum_{S\subset \{ 1, \ldots, m\}} \mes [V_S\cap |f-t|<\kappa].\tag 1.13
$$
Fix $S\subset \{1, \ldots, m\}$.

We consider $f$ on $V_S$ as a function of $(\ve_j)_{j\in S}$ with the other variables $(\ve_j, \ve_j')_{j\in S}$ fixed.
Denoting $g=g(\ve_j; j\in S)\}$ this function on $\{1, -1\}^{|S|}$, we have by our assumption (1.9), for $j\in S$
$$
\align
I_j(g)(\ve_j, j\in S)&= f(\ve_1, \ve_1', \ldots, \ve_{j-1}, \ve_{j-1}', 1, 1, \ve_{j+1},\ve_{j+1}', \ldots,  \ve_n, \ve_n')
\\
&-f(\ve_1, \ve_1', \ldots,\ve_{j-1}, \ve_{j-1}', -1, -1, \ve_{j+1}, \ldots, \ve_n')\\
&\geq\kappa
\endalign
$$
provided
$$\align
(\ve_k)_{k\in S}\in \Omega_j'&= \{(\ve_k)_{k\in S}; \big((\ve_k, \ve_k)_{k\in S}, (\ve_k, \ve_k')_{k\not\in S}) \in\Omega_j\}\\
&= (\Omega_j\cap V_S) \ (\ve_k, \ve_k'; k\not\in S)\subset\{1, -1\}^{|S|}
\endalign
$$
hence only depending on
$ (\ve_k)_{k\in S, k>j}$,  (recall that we fixed the  variables outside $S$).

Applying Lemma 1 to the function $g$ (with $\Omega_j= \{1, -1\}^{|S|}$ for all $j\in S$), we obtain
$$
\align
\# [\omega\in V_S; |f(\omega) -t|<\kappa]&\leq \frac {\# V_S}{|S|^{1/2}}+\sum_{j\in S} \sum_{\ve_k\not= \ve_k', k\not\in S}
\big(2^{|S|} -\# (\Omega_j\cap V_S) (\ve_k, \ve_k'; k\not\in S)\big)\\&= \frac {\# V_S}{|S|^{\frac 12}}+ \sum_{j\in S} \# (V_S\backslash \Omega_j)
.\tag 1.14
\endalign
$$
Summing over $S\subset\{1, \ldots, m\}$ gives
$$
\align
(1.13)&\leq 2^{-m} \sum_{S\subset \{1, \ldots, m\}}\frac 1{|S|^{\frac 12}} +\sum^n_{j=1} \mes(\Omega\backslash\Omega_j)\\
&\lesssim \Big(\frac m2\Big)^{-\frac 12} +n\delta\tag 1.15
\endalign
$$
and hence (1.11).

\bigskip
\noindent
{\bf \S2. Application to the Anderson-Bernoulli model}

Consider the projective action of $SL_2(\Bbb R)$ on $P_1(\Bbb R)\simeq \Bbb T=\Bbb R/\Bbb Z$, defined for $g=\pmatrix a&b\\c&d\endpmatrix
\in SL_2 (\Bbb R)$ by
$$
e^{i\tau_g(\theta)} =\frac{(a \cos \theta +b\sin\theta)+ i(c \cos \theta+d \sin \theta)}{[(a \cos\theta+ b\sin\theta)^2 +
(c \cos \theta +d\sin\theta)^2]^{1/2}}.\tag 2.1
$$
Hence
$$
(\tau_g)'(\theta)=\frac {\sin ^2 \tau_g(\theta)}{(c \cos \theta +d\sin\theta)^2}= \frac 1{[(a\cos \theta + b\sin \theta)^2+ (c \cos\theta + d\sin\theta)^2]^{1/2}}
\tag 2.2
$$
and
$$
\Vert g\Vert^2 \geq (\tau_g)' \geq \frac 1{\Vert g\Vert^2}.\tag 2.3
$$
Consider the Anderson-Bernoulli model (A-B model)
$$
H_\lambda (\ve) =\lambda\ve_n \delta_{nn'}+\Delta\tag 2.4
$$
with $\ve =(\ve_n)_{n\in\Bbb Z} \in \{1, -1\}^{\Bbb Z}$ at small disorder $\lambda>0$ ($\Delta$ stands for the usual lattice Laplacian).

The corresponding transfer operators $M_N(E) \in SL_2(\Bbb R)$ are  given by
$$
\align
M_N=M_N(E;\ve)&=\pmatrix E-\lambda\ve_n& -1\\ 1&0\endpmatrix
\pmatrix E-\lambda\ve_{N-1}&-1\\ 1&0\endpmatrix \cdots \pmatrix E-\lambda\ve_1& -1\\ 1&0\endpmatrix\\
& = \prod^1_{N} g_{_E}(\ve_j).\tag 2.5
\endalign
$$
Considering $\ve_j$ $(1\leq j\leq N)$ as a continuous variable on $[-1, 1]$, and $\partial_j$ the corresponding partial derivative, we
have for the projective action 
$$
\tau_{M_N} = \tau_{g_{_E}(\ve_N)\cdots g_{_E}(\ve_{j+1})^o} \, \tau_{g_{_E} (\ve_j)^o \tau_{g_{_E}}(\ve_{j-1})\cdots g_{_E}(\ve_1)}
$$
and
$$
(\partial_j\tau_{M_N}) (\theta) = \tau'_{g_{_E}(\ve_N)\cdots g_E(\ve_{j+1})} (\tau_{g_E (\ve_j)\cdots g_{_E}(\ve_1)}(\theta)).
(\partial_j\tau_{g_{_E}})(\tau_{g_{_E}(\ve_{j-1})\cdots g_{_E}(\ve_1)}).\tag 2.6
$$
Since
$$
\text{cotg\,} \tau_{g_{_E}(\ve)}(\theta) =(E-\lambda\ve) -\frac {\sin\theta}{\cos \theta}\tag 2.7
$$
we have
$$
(\partial_\ve \tau_{g_{_E}})(\theta)= \lambda.\sin^2 \tau_{g_{_E}(\ve)}(\theta) =\lambda\frac {\cos^2\theta}{\cos^2 \theta+((E-\lambda\ve)\cos\theta-\sin\theta)^2}
\sim \lambda \cos^2 \theta.\tag 2.8
$$
From (2.3), (2.6), (2.8)
$$
\align
(\partial_j\tau_{M_N})(\theta)&\gtrsim \frac \lambda{\Vert g_{_E}(\ve_N)\cdots g_{_E}(\ve_{j+1})\Vert^2}
\cos^2 \tau_{g_{_E}(\ve_{j-1})\cdots g_{_E}(\ve_1)}(\theta)\\
&=\frac \lambda{\Vert M_{N-j} (E;\ve_{j+1}, \ldots, \ve_N)\Vert^2} \cos^2\tau_{M_{j-1} (E, \ve)} (\theta).\tag 2.9
\endalign
$$
In order to deal with the issue of $\cos \tau_{M_{j-1} (E, \ve)}(\theta)$ being small, note that by (2.7), for all $\theta, \ve$
$$
|\cos\theta|+|\cos \tau_{g_{_E}(\ve)}(\theta)|>c.\tag 2.10
$$
Hence (2.9) implies
$$
(\partial_j\tau_{M_N})(\theta)+(\partial_{j+1} \tau_{M_N})(\theta) \gtrsim \frac\lambda{\Vert M_{N-j}(E; \ve_{j+1, \ldots, \ve_N})\Vert^2}\tag 2.11
$$
(for all $\theta$).

In order to fulfill condition (1.9), we need an upperbound on $\Vert M_n (E; \ve_{1}, \ldots, \ve_n)\Vert$.
This function can be analyzed using the Figotin-Pastur  expansion.

Denote
$$
\align
E&= 2 \cos \kappa \qquad (0\leq \kappa\leq\pi)\tag 2.12\\
V_n&= -\frac {\ve_n}{\sin \kappa}\tag 2.13
\endalign
$$
where we assume $\delta_0 <|E|< 2-\delta_0$ and hence $\kappa$ stays away from $0, \frac \pi 2, \pi$ (here $\delta_0$ will be a fixed constant independent of
$\lambda$).

The Figotin-Pastur formula gives
$$
\frac 1N \log\Vert M_N(E, \ve)\Vert =\frac 1{2N} \sum^N_1 \log \big(1+\lambda V_n \sin 2 (\vp_n+\kappa)+ \lambda^2 V_n^2 \sin^2(\vp_n+\kappa)\big)
\tag 2.14
$$
with
$$
\zeta_n =e^{2i\vp_n}\tag 2.15
$$
recursively given by
$$
\zeta_{n+1} =\mu \zeta_n+ i\frac \lambda 2 V_n\frac {(\mu\zeta_n-1)^2}{1-\frac {i\lambda} 2 V_n(\mu\zeta_n-1)}\tag 2.16
$$
and
$$
\mu=e^{2i\kappa}.\tag 2.17
$$
Note that by (2.13), (2.16), $\zeta_n$ only depends on $\ve_{n'}$ for $n'\leq n-1$.

Expanding (2.14), we obtain
$$
\align
(2.14)&=\frac {\lambda^2}{8N} \sum^N_1 V_n^2\tag 2.18\\
& \ +\frac\lambda{2N} \sum_1^N V_n\sin (\vp_n+\kappa)\tag 2.19\\
& \ -\frac {\lambda^2}{4N} \sum_1^N V_n^2 \cos 2(\vp_n+\kappa)\tag 2.20\\
& \ + \frac {\lambda^2}{8N} \sum^N_1 V_n^2 \cos 4 (\vp _n+\kappa)\tag 2.21\\
& \ +0(\lambda^3)
\endalign
$$
and
$$
(2.18) =\frac {\lambda^2}{8\sin^2 \kappa} = \frac{\lambda^2}{2(4-E^2)}.\tag 2.22
$$
By (2.16), (2.17)
$$
\align
&|1-\mu| \ \Big|\sum^N_1 \zeta_n\Big|< 1+0(\lambda N)\\
&\Big|\sum^N_1\xi_n\Big| <\frac{0(\lambda N)}{\sin^2\kappa}\tag 2.23
\endalign
$$
and similarly
$$
\Big|\sum^N_1\zeta_n^2\Big|<\frac{0(\lambda N)}{\sin^2 2\kappa} < 0(\lambda N).\tag 2.24
$$

Writing $\cos 2(\vp_n+\kappa) =\frac 12 (\mu\zeta_n+ \bar\mu\bar\zeta_n), \cos 4 (\vp_n+\kappa) =\frac 12 (\mu^2 \zeta_n^2 +\bar\mu^2 \bar \zeta^2_n)$,
(2.23), (2.24) imply
$$
(2.20), (2.21) =0(\lambda^3).\tag 2.25
$$
Write
$$
\align
(2.19) &=\frac {-\lambda}{2N\sin\kappa}\sum^N_1 \ve_n\sin (\vp_n+\kappa)\\
&=- \frac{\lambda}{2N\sin\kappa}\sum_1^N \ve_n d_n(\ve_{n'}; n'<n)\tag 2.26
\endalign
$$
which is a martingale difference sequence, with
$$
\align
\sum_1^N |d_n|^2& =\sum^N_1\sin^2 2(\vp_n+\kappa)< \frac N2+\frac 12\Big|\sum_1^N \mu^2 \zeta^2_n+\bar\mu^2 \bar\zeta^2_n\Big|\\
&< \Big(\frac 12 +0(\lambda)\Big)N.\tag 2.27
\endalign
$$

In conclusion
$$
\frac 1N\log\Vert M_N(E; \ve)\Vert =\frac {\lambda^2}{8\sin^2\kappa} - \frac \lambda{2N\sin\kappa} \sum_1^N \ve_nd_n+0(\lambda^3)\tag 2.28
$$
and the Lyapounov exponent
$$
L(E) =\frac{\lambda^2} {8\sin ^2\kappa} +0(\lambda^3).\tag 2.29
$$
From martingale theory and (2.28), we get for $a>0$ the large deviation inequality
$$
\align
\mes\Big[\ve\ \, \Big|\frac 1N\log\Vert M_N(E; \ve)\Vert -L(E)\Big|> aL(E)\Big] &< e^{-\big(\frac{a^2\lambda^2}{16\sin ^2\kappa}
+0(\lambda^3)\big) N}\\
&<e^{-\big(\frac {a^2}2 L(E)+0(\lambda^3)\big)N}.\tag 2.30
\endalign
$$
In particular, taking $a>2$ (and $\lambda$ small)
$$
\mes[\ve|\log\Vert M_N (E;\ve)\Vert> aN\lambda^2]< e^{-ca^2 \lambda^2N}.\tag 2.31
$$
Returning to (2.11), take
$$
N\sim \lambda^{-2}.\tag2.32
$$
For $1\leq j\leq N$, it follows from (2.31) that
$$
\align
\mes&[(\ve_{j+1}, \ldots, \ve_N); \Vert M_{N-j} (E; \ve_{j+1}, \ldots, \ve_N)\Vert> e^{C_1(\log N)^{\frac 12}}]\leq\\
&\exp \{[-cC_1^2\lambda^2\frac{\log N}{(\lambda^2(N-j))^2} +0(\lambda^3)] (N-j)\} < e[-c C_1^2 \log N+0(\lambda)] < N^{-C_1}.\tag 2.33
\endalign
$$
Recalling (2.11), we see that Lemma 2 may be applied to the function
$$
f=\tau_{M_N(E; \ve)}(\theta) \text { of } \ve\in \{1, -1\}^N \text { with } \kappa\sim\lambda e^{-2C_1(\log N)^{1/2}}\text { and } \delta< N^{-C_1}
<N^{-10},
$$
for a suitable choice of the constant $C_1$.
\bigskip

Hence, we proved

\proclaim
{Lemma 3} For $\lambda$ small and $N\sim\lambda^{-2}$, we have for fixed $\delta_0<|E|< 2-\delta_0$ and $\theta\in \Bbb T$ arbitrary,
the distributional inequality
$$
\mes[\ve; |\tau_{M_N(\ve)}(\theta) -t|<\lambda e^{-C|\log \lambda|^{\frac 12}}] \leq C\lambda\tag 2.34
$$
for all $t$ (where $C$ is some constant).
\endproclaim

\bigskip

\noindent
{\bf \S3. Dimension of the Furstenberg measure}

Fixing $E$ as above, denote $\nu_E=\nu$ the Furstenberg measure on $\Bbb T$ for the random walk associated with the probability measure on $SL_2(\Bbb R)$
$$
\mu=\frac 12 \, \delta_{\pmatrix E-\lambda &-1\\ 1&0\endpmatrix} + \frac 12  \, \delta_{\pmatrix E+\lambda& -1\\ 1&0\endpmatrix}.\tag 3.1
$$
Thus for all $N$
$$
\int_{SL_2(\Bbb R)} \vp(g) \mu^{(N)} (dg) =\int_{\{1, -1\}^N} \vp\big(M_N(E, \ve)\big) d\ve.\tag 3.2
$$
The measure $\nu$ is $\mu$-stationary i.e.
$$
\nu=\int(\tau_g)_* [\nu] \mu(dg)\tag 3.3
$$
and
$$
\langle\nu, f\rangle =\lim_{N\to\infty} \int f\big(\tau_{M_N(\ve)}(\theta)\big) d\ve\tag 3.4
$$
for all $f\in C(\Bbb T)$ and $\theta\in\Bbb T$.

Our goal is to show that for small $\lambda$, the dimension of $\nu_E$ is close to 1.

The main inequality is the following

\proclaim
{Lemma 4}
Let $h\in SL_2(\Bbb R)$ be arbitrary such that
$$
\Vert h\Vert\sim \lambda^{-\frac 1{10}}.
$$
Let $N\sim\lambda^{-1}$ and $I\subset\Bbb T$ an arbitrary interval of size $|I|<\lambda$. Then
$$
\align
&\int (\tau_{M_N(\ve)h})_* [\nu](I) d\ve\leq\\
&e^{C|\log \lambda|^{\frac 12}} \big\{\max_{|J|<\lambda^{1/10}|I|} \nu(J) +\lambda^{\frac 1{30}} \max_{|J|\leq|I|}\nu(J) +\max_{\lambda^{-\frac
1{10}}<D<\lambda^{-\frac 15}} \frac 1D \max_{|J|<D.|I|} \nu(J)\big\}\tag 3.5
\endalign
$$
with $J$ denoting an interval.
\endproclaim

\noindent
{\it Proof.}
Write
$$
\int\big(M_N(\ve)h\big)_* [\nu] (I) d\ve =\sum_{0\leq k\lesssim N} \int_
{[\Vert M_N(\ve)\Vert\sim 2^k]} \nu(\tau_{h^{-1}} \tau _{M(\ve)^{-1}} (I)\big)
d\ve.\tag 3.6
$$
From (2.31)
$$
\mes [\Vert M_N(\ve)\Vert \sim 2^k]< e^{-ck^2}\tag 3.7
$$
and, if $\Vert M_N(\ve)\Vert\sim 2^k, \tau_{h^{-1}}\tau_{M_N(\ve)^{-1}}(I)$ is contained in an interval $J\in\Bbb T$ of size at most
$\Vert h\Vert^2 4^k|I|$.
Thus the $k^{\text th}$ summands in (3.6) is certainly bounded by
$$
e^{-ck^2} \max_{|J|<4^k\Vert h\Vert^2 |I|} \nu (J).\tag 3.8
$$
Next, restrict $k\lesssim (\log N)^{\frac 12}$ and $\ve$ to $[\Vert M_N(\ve)\Vert\sim 2^k]$.

Let $R_1, \ldots, R_M$ be a partition of $\Bbb T$ in intervals of size $\frac 1M\sim\lambda$.
Estimate
$$
\align
&\int_{[\Vert M_N(\ve)\Vert\sim 2^k]} \nu \big(\tau_{h^{-1}}\tau _{M_N(\ve)^{-1}}(I)\big) d\ve\leq\\
&\sum^M_{m=1} \int_{[\Vert M_N(\ve)\Vert\sim 2^k]} \nu\big(\tau_{h^{-1}}(\tau_{M_N(\ve)^{-1}} (I) \cap R_M)\big) d\ve\\
& \leq \sum^M_{m=1} \mes[\ve; \Vert M_N(\ve) \Vert\sim 2^k \text { and }\tau_{M_N(\ve)}(R_m) \cap I \not=\phi]
\left[ \matrix \max \nu(J)\\ |J|\leq 4^k D_m |I|\endmatrix\right]\tag 3.9
\endalign
$$
denoting
$$
D_m=\max_{\theta\in R_m}|\tau_{h^{-1}}'(\theta)|.\tag 3.10
$$
Fixing some $\theta_m\in R_m$ and $\psi\in I$, $\tau_{M_N(\ve)} (R_m)$ is contained in an $4^k\frac 1M$-neighborhood of $\tau_{M_N(\ve)}(\theta_m)$
and hence
$$
|\tau_{M_N(\ve)}(\theta_m)-\psi|\lesssim \frac {4^k}M+|I|\lesssim \frac{4^k}M\tag 3.11
$$
since $\tau_{M_N(\ve)} (R_m)\cap I\not= \phi$.
In view of Lemma 3
$$
\mes [\ve; (3.11)]< \frac 1M 4^k e^{C(\log N)^{1/2}}
$$
by (2.34) and a suitable partition of the interval $[\psi-\frac{4^k}M, \psi+\frac{4^k}M]$. 
Hence, for $k$ as above
$$
\mes [\ve;\Vert M_N(\ve)\Vert\sim 2^k\text { and }\tau_{M_N(\ve)} (R_m)\cap I\not=\phi] < e^{C(\log N)^{12}}\lambda.\tag 3.12
$$
Let $h^{-1} =\pmatrix a&b\\c&d\endpmatrix$. By (2.2)
$$
\tau_{h^{-1}}'(\theta) =\frac 1{(a\cos\theta+b\sin \theta)^2+(c\cos\theta +d\sin\theta)^2}
$$
and
$$
\frac 1{\Vert h\Vert^2} \lesssim \tau_{h^{-1}}'(\theta) \lesssim \min\Big(\frac 1{\Vert h\Vert^2 \Vert\theta-\theta_h\Vert^2}, \Vert h\Vert^2\Big)\tag
3.13
$$
for some $\theta_h\in\Bbb T$. Thus
$$
\frac 1{\Vert h\Vert^2} \lesssim D_m \lesssim \min \Big[\frac 1{\Vert h\Vert^2 \Vert \theta_m -\theta_h\Vert^2},\Vert h\Vert^2\Big].\tag 3.14
$$
Hence, given $D>0$
$$
\# \{1\leq m\leq M; D_{m} \sim  D\}\lesssim 1+\frac M{\Vert h\Vert D^{1/2}}.\tag 3.15
$$
From (3.12), (3.15), we obtain the following estimate on (3.9)
$$
\align
(3.9)& < e^{C(\log N)^{1/2}} \lambda (\log N)\Big\{ \max_{\Vert h\Vert^{-2}< D<\Vert h\Vert^2} \frac M{\Vert h\Vert D^{1/2}}\big(\max_{|J|< 4^k D|I|}
\nu(J)\big)\Big\}\\
&< e^{C'(\log N)^{1/2}} \big(\max_{\Vert h\Vert^{-2} <D<\Vert h\Vert^2} \frac 1{D^{\frac 12}\Vert h\Vert}\max_{|J|<D|I|} \nu(J)\big)\tag 3.16
\endalign
$$
since $k\lesssim (\log N)^{1/2}$ and writing $J$ as a union of $4^k$ intervals of size at most $D.|I|$.
\bigskip

We distinguish several contributions

(i). For $D<\Vert h\Vert^{-1}$, estimate (3.16) by
$$
e^{C' (\log N)^{\frac 12}} \max_{|J|<\frac{|I|}{\Vert h\Vert}} \nu(J) < e^{C'\vert\log\lambda\vert^{\frac 12}} \max_{|J|<\lambda^{\frac 1{10}} |I|}\nu
(J).\tag 3.17
$$

(ii) For $1>D>\Vert h\Vert^{-1}$, we have $D^{1/2} \Vert h\Vert>\Vert h\Vert^{1/2} \gtrsim \lambda^{-\frac 1{20}}$ and we may bound (3.16) by
$$
\lambda^{\frac 1{30}}\max_{|J|\leq |I|} \nu(J).\tag 3.18
$$

(iii) For $1\leq D\leq \Vert h\Vert$, bound (3.16) by
$$
e^{C' |\log \lambda|^{1/2}} \frac {D^{\frac 12}}{\Vert h\Vert} \max_{|J|\leq|I|}\nu (J) \leq \lambda^{\frac 1{30}} \max_{|J|\leq |I|} \nu (J).\tag 3.19
$$

(iv) For $\Vert h\Vert <D< \Vert h\Vert^2$, estimate by
$$
\frac {e^{C |\log \lambda|^{1/2}}} D \max_{|J|<D|I|} \nu(J).\tag 3.20
$$
Collecting the contributions (3.17) - (3.20) gives (3.5).
This proves Lemma 4.
\bigskip

Next, returning to (3.3), writing $\mu =\frac 12 \delta_{g_1} +\frac 12\delta_{g_2}$
we make the following construction. Assume
$$
\nu =\int (\tau_g)_* [\nu] \mu_1(dg)\tag 3.21
$$
where $\mu_1$ is some discrete probability measure $SL_2(\Bbb R)$, such that
$$
\Vert g\Vert < 2\lambda^{-\frac 1{10}} \text { for } g\in \supp\mu_1.\tag 3.22
$$
If $g\in\supp \mu_1$ and $\Vert g\Vert <\lambda^{-\frac 1{10}}$, write by (3.3)
$$
(\tau_g)_* [\nu] =\frac 12 (\tau_{g g_1})_* [\nu] +\frac 12 (\tau_{g g_2})_* [\nu].
$$
Define then
$$
\mu_2 =\sum_{\Vert g\Vert\geq \lambda^{-\frac 1{10}}}\mu_1(g) \delta_g +\frac 12 \sum_{\Vert g\Vert<\lambda^{-\frac 1{10}}} \mu_1(g) (\delta_{g g_1}+\delta{ g
g_2})
$$
still satisfying (3.21).

From the positivity of the Lyapounov exponent, an iteration of this process will clearly produce a discrete probability measure $\tilde\mu$ on
$SL_2(\Bbb R)$ s.t.
$$
\nu = \int (\tau_g)_* [\nu] \tilde\mu (dg)\tag 3.23
$$
and
$$
\lambda^{- \frac  1{10}} <\Vert g\Vert < 2\lambda^{-\frac 1{10}} \text { for } g\in \supp \tilde\mu.\tag 3.24
$$
Taking $N\sim \lambda^{-2}$, and since also by (3.2)
$$
\nu = \int (\tau_{M_N(\ve)})_* [\nu] d\ve
$$
(3.23) implies
$$
\nu =\int\Big[\int (\tau_{M_N(\ve)h})_* [\nu] d\ve\Big]\tilde\mu (dh).\tag 3.25
$$

From(3.25) and Lemma 4, we conclude the following inequality.

\proclaim
{Lemma 5} For $I\subset\Bbb T$ an interval of size at most $\lambda$, we have
$$
\nu(I) \leq e^{C|\log \lambda|^{\frac 12}} \Big[\max_{|J|<\lambda^{\frac 1{10}}|I|} \nu(J) +\max_{\lambda^{-\frac 1{10}} <D<\lambda^{-\frac 15}}
\frac 1D \max_{|J|<D|I|} \nu(J) \Big].\tag 3.26
$$
\endproclaim

If we iterate (3.26) $r$-times, assuming $\lambda^{-\frac r5}|I|<\lambda$, we obtain
$$
\nu(I) \leq 2^r e^{C\vert\log \lambda|^{\frac 12} r} \frac 1{D_1} \nu (J)\tag 3.27
$$
for some interval $J$ of size $|J|<D_1 \delta_1|I|$ where $D_1> 1, 0<\delta_1< 1$ and $D_1\delta_1^{-1}<\lambda^{-\frac r{10}}$.

From random matrix product theory it is known that the Furstenberg measure $\nu$ has some positive dimension $\alpha>0$.
Hence the right side of (3.27) is at most
$$
\lesssim C^{|\log\lambda|^{\frac 12} r} \frac 1{D_1} |J|^\alpha < C^{|\log\lambda |^{\frac 12} r} \frac {\delta_1^\alpha}{D_1^{1-\alpha}}
|I|^\alpha < 
C^{|\log \lambda|^{\frac 12} r} \lambda^{\frac r{20} \ \min (\alpha, 1-\alpha)} |I|^\alpha.\tag 3.28
$$
Assuming $\gamma>0$ some constant (independent of $\lambda=o(1)$) satisfying
$$
\gamma<\alpha< 1-\gamma\tag 3.29
$$
(3.28) and the restriction on $r$ would imply for $\lambda<\lambda(\gamma)$
$$
\nu(I)<(C^{|\log \lambda|^{\frac 12} }\lambda^{\frac \gamma{20}})^r |I|^\alpha < \lambda^{\frac {\gamma r}{30}} |I|^\alpha
$$
and
$$
\nu (I) <\Big(\frac{|I|}\lambda\Big)^{\frac \gamma 6} |I|^\alpha \lesssim |I|^{\alpha+\frac \gamma 6}.\tag 3.30
$$
But (3.30) would give that $\nu$ has dimension at least $\alpha+\frac \gamma 6$, a contradiction.

Thus in order to prove

\proclaim
{Theorem 1} Assuming $\delta_0<|E|< 2-\delta_0$, the dimension of the Furstenberg measure $\nu_E^{(\lambda)}$ for the A-B model is at least
$\alpha(\lambda)\overset {\lambda\to 0}\to\rightarrow 1$.
\endproclaim

It will suffice to have a uniform lower bound in $\lambda$ for dim\,$\nu_E^{(\lambda)}$.

This is what we establish next.

\proclaim
{Lemma 6} Under the assumption of Theorem 1, $\dim \nu_E^{(\lambda)} > \gamma>0$ with $\gamma$ independent of $\lambda$.
\endproclaim

\noindent
{\it Proof.} We make the following observation.
Write $M=M_N(E; \ve)$ as
$$
M = \frac {(v_-^\bot \otimes v_+) \lambda_++ (v_+^\bot \otimes v_-)\lambda_+^{-1}}{\langle v_+, v_-^\bot \rangle}\tag 3.31
$$
with $v_+$ (resp. $v_-$) the expanding (resp. contracting) direction. Hence
$$
\Vert M\Vert\sim \frac {|\lambda_+|}{|v_+\wedge v_-|}.\tag 3.32
$$

For unit vectors $u, w\in\Bbb R^2$, we deduce from (3.31) that
$$
\align
\frac {\Vert M u\Vert}{\Vert M\Vert}&= \Vert\langle v_-^\bot, u\rangle v_+ +\lambda^{-2}_+ \langle v_+^\bot, u\rangle v_-\Vert = \\
(1+\lambda_+^{-2})|\langle v^\bot_-, u\rangle|+ 0\Big(\frac{|v_+\wedge v_-|} {\lambda^2_+}\Big) &\overset{(3.32)}\to \leq
(1+\lambda_+^{-2})|\langle v_-^\bot, u\rangle|+0\Big(\frac 1{\Vert M\Vert}\Big)\tag 3.33
\endalign
$$
and
$$
\align
\frac {|\langle M u, w\rangle|}{\Vert M\Vert} &= |\langle v_-^\bot, u\rangle \langle v_+, w\rangle
+\lambda^{-2} \langle v_+^\bot, u\rangle \langle v_-, w\rangle|\\
&\geq (1+\lambda^{-2}) |\langle v_-^\bot, u\rangle| \ |\langle v_+, w\rangle| -2 \lambda^{-2}_+|v_+\wedge v_-|\\
&\geq |\langle v_-^\bot, u\rangle| \ |\langle v_+, w\rangle| +0\Big(\frac 1{\Vert M\Vert}\Big).\tag 3.34
\endalign
$$
Hence, given an arc $I$ of size $\eta$ centered at $\nu$
$$
\align
&\Bbb P\Big[\ve; v_-\in I \text {  where $v_-$ is contracting direction  of $M_N(\ve)\Big]\leq ^{(3.33)}$}\\
&\Bbb P\Big[\ve; \frac {\Vert M_N(\ve) u\Vert}{\Vert M_N(\ve)\Vert} < 2\eta+0\Big(\frac 1{\Vert M_N(\ve)\Vert}\Big)\Big]\leq\\
&\Bbb P\big[\ve; \Vert M_N(\ve)\Vert < e^{\frac {\lambda^2}{20}N}\Big] +\Bbb P\Big[\ve; \frac{\Vert M_N(\ve) u\Vert}{\Vert M_N(\ve)\Vert}
<3\eta\Big]=\\
& (3.35)+ (3.36)
\endalign
$$
provided
$$
\eta> e^{-\frac{\lambda^2}{20} N}.\tag 3.37
$$
Recalling (2.29), (2.30), we have
$$
\mes \Big[ \ve\ \Big| \frac 1N \frac {\log\Vert M_N(\ve)\Vert}{L(E)} -1\Big|> a\Big] < e^{\big(-\frac{a^2}2 L(E)+0(\lambda^3)\big)N}
\tag 3.38
$$
with $\frac{\lambda^2}8< L(E)<0(\lambda^2)$.

Hence,
$$
(3.35)< e^{-(\frac 1{50}\lambda^2 +0(\lambda^3))N} < e^{-\frac 1{60}\lambda^2 N} \overset{(3.37)}\to < \eta^{\frac 13}\tag 3.39
$$
for $\lambda$ small enough.

Next, we point out that in the analysis (2.14)-(2.28), the formula (2.28) is equally valid for $\frac 1N\log\Vert M_N(E; \ve)(u)\Vert$,
with $u\in S^1$ arbitrary (as a consequence of the argument).  Thus we can write
$$
\frac 1N\log\Vert M_N(\ve)\Vert =L(E) -\frac\lambda{2N\sin\kappa} \sum^N_1 \ve_n d_n+ 0(\lambda^3)\tag 3.40
$$
and
$$
\frac 1N\log \Vert M_N(\ve) (u)\Vert =L(E) -\frac\lambda{2N\sin\kappa}\sum^N_1 \ve_n d_n' +0(\lambda^3)\tag 3.41
$$
so that
$$
\log \frac {\Vert M_N(\ve)\Vert}{\Vert M_N(\ve) (u)\Vert} =\frac \lambda{2\sin \kappa} \sum_1^N \ve_n (d_n'-d_n)+0(N\lambda^3)\tag 3.42
$$
where $d_n, d_n'$ depend on $\ve_1, \ldots , \ve_{n-1}$.

Letting $1>t> 0$ be a parameter, write 
$$
\align
(3.36) &<(3\eta)^t \int\Big(\frac{\Vert M_N(\ve)\Vert}{\Vert M_N(\ve) (u)\Vert}\Big)^t d\ve< (3\eta)^t e^{0(N\lambda^3 t)}
\int e^{\frac {\lambda t}{2\sin \kappa} \sum_1^N \ve_n (d_n'-d_n)} d\ve
\\
&{}\\
&<(3\eta)^t e^{0(N\lambda^3t)} e^{C\lambda^2t^2N}\tag 3.43
\endalign
$$
where the constant $C$ only depends on $E$.

Choosing $N$ s.t.
$$
\eta\sim e^{-\frac {\lambda^2}{10^3}N}\tag 3.44
$$
we satisfy (3.37), and it follows from (3.43) and appropriate choice of $t$, that
$$
(3.36)<(3\eta)^{t-C(\lambda t+ t^2)}< \eta^{c_1}
$$
(again for $\lambda$ small enough) and with $c_1>0$ independent of $\lambda$.

Hence, we showed that with $N$ satisfying (3.44)
$$
\mes[\ve; v_-\in I\text { where $v_-$ is contracting vector of $M_N(\ve)]<\eta^{c_1}$}.\tag 3.45
$$
Since $v_+$ is the contracting vector of $M_N(\ve)^{-1}$, we obtain a similar statement for the expanding vector.
Therefore, given any pair of $\eta$-intervals $I_+, I_-$ in $S^1$, we proved that
$$
\align
&\mes [\ve; v_+\in I_+, v_- \in  I_- \\
&\text{ with $v_+ $ (resp. $v_-)$ expanding (resp. contracting) direction of }
M_N(\ve)] < 2\eta^{c_1}\tag 3.46
\endalign
$$
for $N$ satisfying (3.44).

Returning to (3.34), we have
$$
\align
&\Bbb P[\ve; \frac {|\langle M_N(\ve) u, w\rangle|}{\Vert M_N(\ve)\Vert}< \eta_1]\leq\\
& \Bbb P[\ve; \Vert M_N(\ve)\Vert< 1/{\eta_1}]+ \Bbb P[\ve; |\langle v_-^\bot, u\rangle|\lesssim \sqrt{\eta_1}]+\Bbb P[\ve; |\langle v_+,
w\rangle|\lesssim \sqrt{\eta_1}].\tag 3.47
\endalign
$$
Taking $\eta =\eta_1^{\frac 12}$ and $N$ as in (3.44), the last 2 terms in (3.37) are at most $0(\eta_1^{\frac 12 c_1})$ by 
(3.46), while the first term is bounded by $\mes[\ve ; \Vert M_N(\ve)\Vert < e^{\frac 1{500\lambda^2 N}}] < 
e^{-\frac 1{60} \lambda^2 N} <\eta_1$ by
(3.38).

Hence
$$
(3.47)\lesssim \eta_1^{\frac 12c_1} \text { with $\eta_1 \sim e^{-\frac{\lambda^2}{500} N}$}.\tag 3.48
$$
Returning to the Furstenberg measure $\nu =\nu_E^{(\lambda)}$, we have for $I\subset\Bbb T$ a small arc of size $\eta_1$, by (3.4)
$$
\nu(I) =\lim_{N'\to\infty} \Bbb P\Big[\ve| \frac{M_{N'}(\ve) e_1}{\Vert M_{N'}(\ve) e_1\Vert} \in I\Big].
$$

Take $N$ as in (3.48) and $N'>N$. If $w$ denotes the center of $I$, then
$$
\frac {|\langle M_{N'} e_1, w^\bot\rangle|}{\Vert M_{N'} e_1\Vert} <\eta_1.\tag 3.49
$$
Fix $\ve_1, \ldots, \ve_{N'-N} $ and let $u=\frac{M_{N'-N}(\ve_1, \ldots, \ve_{N'-N})(e_1)}{\Vert M_{N'-N_1}, \ldots,
(\ve_{N'-N})(e_1)\Vert}$.
We have
$$
\frac{M_{N'}e_1}{\Vert M_{N'}e_1\Vert} =\frac{M_N(\ve_{N'-N+1}, \ldots, \ve_{N'})(u)}
{\Vert M_N(\ve_{N'- N+1}, \ldots, \ve_{N'})u\Vert}.
$$
Thus (3.49) implies
$$
\frac{|\langle M_N(\cdots) u, w^\bot\rangle|}{\Vert M_N(\cdots)\Vert} <\eta_1\tag 3.50
$$
for which the measure in $\ve_{_{N'-N+1}}, \ldots, \ve_{N'}$ is at most $\eta_1^{\frac 12 c_1}$ by(3.48).
Therefore
$$
\nu(I)\lesssim |I|^{\frac 12 c_1}.\tag 3.51
$$
This proves that $\dim \nu\geq \frac 12 c_1$, uniformly in $\lambda$.
Hence we establish Lemma 6.

This also completes the proof of Theorem 1.

\bigskip

\noindent
{\bf \S4. Density of states}

Let $u, w\in S^1$, $\eta>0$ small.
It follows from (3.34) that
$$
\align
&\lim_{N\to\infty} \mes \Big[ ve; \frac{|\langle M_N(\ve)u, w\rangle|}{\Vert M_N(\ve)\Vert} <\eta\Big] \leq\\
&\lim_{N\to\infty} \mes \Big[\ve; |\langle v_+, w\rangle| .|\langle v_-^\bot, u\rangle |< \eta 
\text { with $v_+, v_-$ the eigenvectors of $M_N(\ve)\Big]=$}\\
& \lim_{N\to\infty} \mes \Big[(\ve, \ve'); |\langle v_+, w\rangle|. |\langle v_+', u^\bot\rangle|<\eta
\text { with $v_+$ (resp $v_+'$) expanding direction of }\\
&\text{$M_N(\ve)$, (resp $M_N(\ve'))\Big]$}\\
&\lesssim \log\frac 1\eta .\max_{\eta_1. \eta_2 =\eta} \nu_E(I_{\eta_1} (w^\bot)).\nu_E (I_{\eta_2} (u))\ll \eta^\gamma.\tag 4.1
\endalign
$$
where we used (3.4) and the independence of $v_+, v_-$ for $N\to\infty$ as functions of $\ve$.

Here $\gamma<\dim \nu_E^{(\lambda)}$ and $\gamma =\gamma(\lambda)\to 1$ for $\lambda\to 0$.

It is easily seen that (4.1) implies that for given $K>1$ and taking $N$ large enough (depending on $K$)
$$
\max_{u, w\in S^1} \Bbb E\Big[\frac {\Vert M_N\Vert}{|\langle M_N u, w\rangle |}\wedge K\Big] < K^{1-\gamma}.\tag 4.2
$$
Here $M_N=M_N(E)$ and (4.2) remains clearly valid replacing $E$ by $z= E+iy$ with $0< y< y_N$ small enough (depending on $N$) and taking for $u, w$
unit vectors in $\Bbb C^2$.

Next, take $N' >N$ and consider
$$
\frac {\Vert M_{[0, N']} (z; \ve)\Vert. \Vert M_{]1N', 2N']} (z; \ve)\Vert}
{\Vert M_{[0, 2N']}(z, \ve)\Vert}.\tag 4.3
$$
Fixing $\ve_{N'+1}, \ldots, \ve_{2N'}$, we obtain a unit vector $\zeta\in\Bbb C^2$ (depending on these variables) such that
$$
(4.3) =\frac {\Vert M_{[0, N']}(z, \ve)\Vert}{\Vert M_{[0, N']} (z, \ve)(\zeta)\Vert}\tag 4.4
$$
and
$$
(4.4) \lesssim \sum_{i, j=1, 2} \frac {|\langle M_{[0, N']}(z, \ve) e_i, e_j\rangle|}{|\langle M_{[0, N'']}(z, \ve)\zeta, e_j\rangle|}.
$$
Fix also $\ve_1, \ldots, \ve_{N'-N}$ and let $\zeta_1$ be a unit vector in $\Bbb C^2$ with $\zeta_1$ parallel to $ M^*_{[0, N'-N]} (z, \ve) e_j$.
Hence
$$
\frac{|\langle M_{[0, N']}(z, \ve) e_i, e_j\rangle |}{|\langle M_{[0, N']}(z, \ve)\zeta, e_j\rangle|}\leq
\frac{\Vert M_{[N'-N+1, N']} (z, \ve)\Vert}{|\langle M_{[N'-N+1, N']}(z, \ve)\zeta, \zeta_1\rangle|}\tag 4.5
$$
where the vectors $\zeta, \zeta_1$ do not depend on $\ve_{N'-N+1}, \ldots, \ve_{N'}$.

Thus
$$
\min \big((4.3), K\big) \lesssim \min \big((4.5), K\big).\tag 4.6
$$
Taking expectation of (4.6) in $\ve_{N'-N+1}, \ldots, \ve_{N'}$ (with other variables fixed), (4.2) and subsequent remark,
give an estimate
$$
\Bbb E_{\ve_{N'- N+1}, \ldots, \ve_{N'}} [(4.6)] \lesssim K^{1-\gamma}.\tag 4.7
$$
Hence, also
$$
\Bbb E[\min \big((4.3), K\big)]\lesssim K^{1-\gamma}\tag 4.8
$$
valid for $z =E+iy$ with $y>0$ small enough (depending on $K$) and $N'>N'(K)$.

Denoting $\Cal N$ the IDS, recall that
$$
\bar\partial\Cal N(z) =\Bbb E[G(0, 0, z)] \qquad z= E + iy
$$
where $G(z) =(H-z)^{-1} $ is the Green's function and $\Cal N(z)$ the harmonic extension of $\Cal N$ to Im$\,z>0$.

Fix $z, \text{Im} z>0$. 
Then from the resolvent identity and positivity of the Lyapounov exponent, we obtain
$$
G(0, 0, z) =\lim_{\Sb \Lambda =[-a, b]\\a, b\to\infty\endSb} G_\Lambda (0, 0, z) \ \text { a.s}
$$
and, by Cramer's rule
$$
|G(0, 0, z)|\leq{\underline{\lim}}_{N'\to\infty} \frac {\Vert M_{[-N', 0]} (z, \ve)\Vert \ \Vert M_{[0, N']}(z, \ve)\Vert}
{\Vert M_{[-N', N']} (z; \ve)\Vert}.\tag 4.9
$$
Hence by (4.8)
$$
\Bbb E[|G(0, 0, z)|\wedge K]\leq {\underline\lim}_{N'\to \infty} \Bbb E[\cdots \wedge K ]\lesssim K^{1-\gamma}\tag 4.10
$$
if $y>0$ is small enough (depending on $K$). Letting $y\to 0$ we get
$$
\Bbb E[|G(0, 0, E+ io)|\wedge K]<K^{1-\gamma}.\tag 4.11
$$
It follows from (4.11) that for $0<\gamma_1<\gamma$
$$
\Bbb E[|G(0, 0, E+io)|^{\gamma_1}]R<C.\tag 4.12
$$
Recall that we assumed $\delta_0 < |E|< 2-\delta_0$.
Using the subharmonicity of $|G(0, 0, z)|^{\gamma_1}$ on Im\,$z>0$, we deduce from (4.12) that for fixed $z= E+iy, y>0$
$$
|\bar \partial \Cal N(z)|\leq \Bbb E[|G(0, 0, z)|]<\frac 1{y^{1-\gamma_1}} \Bbb E[|G(0, 0, z)|^{\gamma_1}]<\frac C{y^{1-\gamma_1}}.\tag 4.13
$$

Hence $\Cal N$ is $\gamma_1$-H\"older for all $\gamma_1<\gamma$.

This proves 

\proclaim {Theorem 2} For $\delta_0<|E|< 2-\delta_0$, the IDS of the  A-B model with $\lambda$-disorder is $s$-H\"older regular, with
$s\to 1$ for $\lambda\to 0$.
\endproclaim

\bigskip

\noindent
{\bf \S5. Further comments}

If one aims at going further and prove the Lipschitz regularity of the IDS, it seems reasonable to prove that the Furstenberg measures on the projective
line $P_1(\Bbb R)\simeq \Bbb T$ are at least absolutely continuous.  This is far from an obvious issue.
In fact, it was conjectured in [K-L] that if $\nu$ is a finitely supported probability measure on $SL_2(\Bbb R)$, then its Furstenberg measure on $P_1(\Bbb
R)$ is always singular.
This conjecture was disproved in [BPS] using a probabilistic construction reminiscent of random Bernoulli-convolutions.
An explicit example was given recently in [B2], based on a construction from [B3] \big(that relies on an extension of the spectral gap theory from [BG1] for
$SU(2)$ to $SL_2(\Bbb R)$\big).
A rough description is as follows.
One produces a finite subset $\Cal G\subset SL_2(\Bbb R) \cap \Mat_{2\times 2}(q)$, $q$ a fixed large integer, such that $\log (\#\Cal G)\sim \log q$,
$\Cal G$ generates freely the free group on $\#\Cal G$ generators and moreover $\Cal G$ is contained in a small neighborhood of the 
identity (depending on $q$).
Denoting
$$
\nu =\frac 1{(\# \Cal G)} \sum_{g\in\Cal G} \delta_g\tag 5.1
$$
the probability measure on $SL_2(\Bbb R)$, it is shown that there is a spectral gap for the projective representation $\rho$, in the following sense.
Let $f\in L^2(\Bbb T)$, $\Vert f\Vert_2=1$ and assume $\hat f(n)=0$ for $|n|>K$, where $K=K(q)$ is a sufficiently large constant.
Then
$$
\frac 1{(\# \Cal G)} \Big\Vert\sum_{g\in \Cal G} \rho_g f\Big\Vert_2< \frac 12\tag 5.2
$$
where $\rho_g f=(\tau_g')^{-\frac 12} (f\circ\tau_g)$ and $\tau_g$ the action on $\Bbb T$ defined for $g=\pmatrix a&b\\ c&d\endpmatrix$ by
$$
e^{i\tau_g(\theta)} =\frac {(a\cos \theta +b\sin\theta)+i(c \cos \theta+ d\sin\theta)}{[(a\cos \theta+ b\sin\theta)^2+(c\cos \theta +
d\sin\theta)^2]^{\frac 12}}.\tag 5.3
$$

Since $g\in\Cal G$ are close to identity, (5.2) clearly implies that for $f$ as above
$$
\frac 1{(\# \Cal G)}\Big\Vert\sum _{g\in\Cal G} (f\circ \tau_g)\Big\Vert_2 <\frac 34.\tag 5.4
$$
From (5.4), one may then derive easily that $\nu$ has an a.c.
Furstenberg measure with $C^k$-density, where $k$ can be made arbitrarily large.

It should be pointed out that the contractive properties (5.2), (5.4) are not exploiting hyperbolicity (at least in the usual sense), as the Lyapounov
exponent of the random matrix product corresponding to $\nu$ is small.

Returning to the A-B-model with small $\lambda$, denote
$$
\nu_{\lambda, E} =\frac 12 \delta_{\pmatrix E+\lambda&-1\\ 1&0\endpmatrix} +\frac 12 \delta_{\pmatrix E-\lambda &-1\\ 1&0\endpmatrix}\tag 5.5
$$
and $\nu^{(\ell)}_{\lambda, E}$ its $\ell$-fold convolution.
It seems reasonable to believe that
$$
\Big\Vert\sum_g \nu^{(\ell)}_{\lambda, E} (g) (f\circ \tau_g)\Big\Vert_2 \leq \frac 12\Vert f\Vert_2\tag 5.6
$$
for $f\in L^2(\Bbb T), \hat f(n) =0$ for $|n|>K(\lambda)$ and where $\ell$ is some positive integer independent of $\lambda$, or at least $\ell
=o(\lambda^{-2})$.
Such property would then again imply a.c. and a certain smoothness of the Furstenberg measure.
Unfortunately, available technology to establish spectral gaps (as developed in [BG1]) so far requires algebraic matrix elements of bounded height and hence does
not apply to (5.5).

One may however combine the methods from [BG1] with those of [S-T] to prove the following result, which seems new (compare also with the results
from [K-S]).

\proclaim
{Theorem 3} Consider a random Schr\"odinger operator $H=\Delta+V$ on $\Bbb Z$ where $V=(V_n)_{n\in\Bbb Z}$ are $i.i.d$'s with distribution given by
a compactly supported measure $\beta$ on $\Bbb R$ of positive dimension.
Thus there is $\kappa>0$ s.t.
$$
\beta(I) \lesssim |I|^\kappa \text { for $I\subset\Bbb R$ an interval}.\tag 5.7
$$
Then $H$ has $C^\infty$ density of states.
\endproclaim

We sketch the argument.

For fixed $E$, let $\mu_E$ be the probability measure on $SL_2(\Bbb R)$ obtained as image measure of $\beta$ under the map
$$
v\mapsto \pmatrix E-v&-1\\ 1&0\endpmatrix.\tag 5.8
$$
Following [S-T], it will suffice to show that, for some fixed convolution power $\ell$, the measure $\mu_1=\mu_E^{(\ell)}$ on $SL_2(\Bbb R)$ gives
a smoothing convolution operator on $P_1(\Bbb R)$.
Thus there is some $\alpha>0$ s.t. for $f\in H^s(\Bbb T), s\geq 0$
$$
\Big\Vert \int (f\circ \tau_g)\mu_1(dg)\Big\Vert_{H^{s+\alpha}} \lesssim \Vert f\Vert_{H^s}\tag 5.9
$$
(where $H'$ denotes the usual Sobolev space with norm $\Vert f\Vert_{H^s}=\big(\sum(1+|n|)^{2s} |\hat f(n)|^2\big)^{\frac 12}$\big).

Denoting $x=E-v$, one has
$$
\pmatrix x&-1\\ 1&0\endpmatrix \pmatrix y&-1\\ 1&0\endpmatrix \pmatrix z&-1\\1&0\endpmatrix =\\
\pmatrix xyz-x-z & 1-xy\\ yz-1 &-y\endpmatrix\tag 5.10
$$
and recalling (5.7), one sees that $\mu_E^{(3)}$ certainly has the property that
$$
\mu_E^{(3)} (\frak S_\delta)\lesssim \delta^{\kappa'} \text { for all $\delta>0$}\tag 5.11
$$
if $\frak S$ is a proper algebraic subvariety of $SL_2(\Bbb R)$ of bounded degree and $\frak S_\delta$ denotes a $\delta$-neighborhood of $\Cal G$.
Here $\kappa' >0$ depends on the degree bound.

Let $P_\delta$, $\delta>0$, denote an approximate identity on $SL_2(\Bbb R)$.
Using (5.11), an extension of the `flattening Lemma' from [BG1] to $SL_2(\Bbb R)$ (note that, up to complexification, $SU(2)$ and
$SL_2(\Bbb R)$ have the same Lie-algebra and our analysis is local), permits us to conclude the following.

\proclaim
{Lemma 7} Fix some $0<\ve<1$.
There is $\ell =\ell (\ve) \in\Bbb Z_+$ s.t. for all $\delta>0$, we have
$$
\Vert\mu_E^{(3\ell)}*P_\delta\Vert_\infty <\delta^{-\ve}\tag 5.12
$$
(in particular, $\mu_E^{(3\ell)}$ has dimension at least $3-\ve$).
\endproclaim

This is the crucial step, depending on `arithmetic combinatorics' in groups (see [BG1] and related refs for more details).

Taking $\ve=10^{-3}$ and $\ell =\ell (\ve)$ given by Lemma 7, we can now prove that $\mu_1= \mu_E^{(3\ell)}$ satisfies (5.9).
This will clearly be a consequence of the following statement.

\proclaim
{Lemma 8} Let $f\in L^2(\Bbb T), \Vert f\Vert _2 =1$ and supp\,$\hat f\subset [2^k, 2^{k+1}]\cup [-2^{k+1}, -2^k]$ with $k$ sufficiently large.
Then
$$
\Big\Vert \int (f\circ \tau_g)\mu_1(dg)\big\Vert _2< 2^{-k\kappa}.\tag 5.13
$$
for some $\kappa>0$.
\endproclaim

\bigskip
\noindent
{\bf Proof of Lemma 8.}

We summarize the argument from [B2].

Denote $G=SL_2(\Bbb R)$ and take $\delta= 4^{-k}$, so that, by assumption on $f$, we may replace the left side of (5.13) by
$$
\Big\Vert \int (f\circ \tau_g)(\mu_1*P_\delta)(dg)\Big\Vert_2.\tag 5.14
$$
Using (5.12), one gets
$$
(5.14)^2\lesssim \delta^{-2\ve} \iint_{G\times G}|\langle f\circ\tau_{g_1}, f\circ \tau_{g_2}\rangle|\Omega(g_1)\Omega(g_2) dg_1 dg_2
$$
with $0\leq \Omega\leq 1$ a suitable compactly supported function on $G$ (depending on the support of $\beta$).
Next, by Cauchy-Schwarz
$$
(5.14)^4 \lesssim \delta^{-4\ve}\iiiint_{G\times G\times \Bbb T\times \Bbb T} f(\tau_{g_1} x)\bar f(\tau_{g_2} x) f(\tau_{g_1}y)
f (\tau_{g_2}y) \Omega(g_1) \Omega(g_2) dg_1
dg_2 dxdy.\tag 5.15
$$
To estimate (5.15), proceed as follows. Fix $x, y \in\Bbb T$ and $g_1 \in G$ and consider the integral in $g_2$
$$
\int \bar f(\tau_g x) f(\tau_g y)\Omega(g) dg.\tag 5.16
$$
The point here is that if one specifies $\tau_g x\in\Bbb T$, there remains an average in $\tau_g y$ to be exploited, when integrating in $g$ (unless $x$ and $y$
are very close).  More precisely, if $\Vert x-y\Vert < 2^{-h/10}$, then
$$
|(5.16)|< 2^{-k} \Vert f\Vert_1^2
$$
and the contribution in (5.15) is at most
$$
\leq \delta^{-4\ve} 2^{-k} \Vert f\Vert_1^4 < 2^{-k/2}.
$$
The contribution of $\Vert x-y\Vert < 2^{-k/10} $ in (5.15) is easily estimated by
$$
\align
&\iint_{\Vert x-y\Vert< 2^{-k/10}} \Big[\int_G |f(\tau_g x)| \ |f(\tau_g y)|\Omega(g) dg\Big]^2 dxdy\leq\\
&\iint_{\Vert x-y\Vert < 2^{-k/10}} \Big[\int_G|f(\tau_g x)|^2 \Omega (g) dg\Big] \Big[\int_G |f(\tau_g y)|^2 \Omega(g)dg\Big]dxdy\\
&{}\\
&\lesssim 2^{-k/10}\Vert f\Vert^4_2
\endalign
$$
and (5.13) follows.

\enddocument


\Refs

\widestnumber\no{XXXXXXXXXX}
\ref
\no{[B1]} \by J.~Bourgain
\paper On localization for lattice Schr\"odinger operators involving Bernoulli variables
\jour LNM, 1850, 77--99 (2004)
\endref

\ref
\no {[B2]} \by J.~Bourgain
\paper Finitely supported measures on $SL_2(\Bbb R)$ which are absolutely continuous at infinity
\jour preprint 2010
\endref

\ref
\no{[B3]} \by J.~Bourgain
\paper Expanders and dimensional expansion
\jour C.R. Math.~Acad.~Sci. Paris 347 (2009), no 7--8, 357--362
\endref

\ref\no{[BG1]} \by J.~Bourgain, A.~Gamburd
\paper On the spectral gap for finitely-generated subgroups of $SU(2)$
\jour Inv.~Math. 171 (2008), no 1, 83--121
\endref

\ref\no[{BG1]}\by J.~Bourgain, A.~Gamburd
\paper Spectral gap in $SU(d)$
\jour C.R.~Math.~Acad.~Sci. Paris 348 (2010), no 11--12, 609--611
\endref

\ref\no {[B-K]} \by J.~Bourgain, C.~Kenig
\paper
On localization in the continuous Anderson-Bernoulli model in higher dimension
\jour Inv.~Math. 161 (2005), no 2, 389--426
\endref

\ref\no{[B-P-S]}\by B.~Barany, M.~Pollicott, K. Simon
\paper Stationary measures for projective transformations: the Blackwell and Furstenberg measures 
\jour preprint 2010
\endref

\ref\no {[C-K-M]}\by R.~Carmona, A.~Klein, F.~Martinelli
\paper Anderson localization for Bernoulli and other singular potentials
\jour CMP 108, 41--66 (1987)
\endref

\ref\no {[K-L]}\by V.~Kaimanovich, V.~Le Prince
\paper Matrix random products with singular harmonic measure
\jour Geom.~Ded. (2010)
\endref

\ref\no {[K-S]} \by A.~Klein, A.~Speis
\paper Regularity of the invariant measure and the density of states in the one-dimensional Anderson model
\jour JFA 88 (1990), no 1, 211--227
\endref

\ref\no {[S-T]} \by B.~Simon, M.~Taylor
\paper Harmonic analysis on $SL_2(\Bbb R)$ and smoothness of the density of states in the one-dimensional Anderson model
\jour CMP, 101, no 1 (1985), 1--10
\endref

\ref\no{[S-V-W]}\by C.~Shubin, T.~Vakilian, T.~Wolff
\paper Some harmonic analysis questions suggested by Anderson-Bernoulli models
\jour GAFA 8 (1988), 932--964
\endref
\endRefs
\enddocument